\documentclass[12pt,leqno]{article}
 
\def\w{\dot{w}}
\def\v{{\rm v}} 
 
\def\int{\mathbb{Z}}

\def\Ue{{\cal U}_{\varepsilon}({\mathfrak g})}
\def\g{\mathfrak g}

\newcommand{\Od}{\mathcal{O}}
\newcommand{\NN}{\mathcal{N}}
\newcommand{\Aut}{\mathrm{Aut}}
\newcommand{\lf}{\mathfrak{l}}

\newcommand{\so}{\mathfrak{so}}
 
\def\proof{{\bf Proof. }}

\def\pf{\proof}

\def\mathM{{\mathcal M}}
\def\mathT{{\mathcal T}}

\def\OO{{\cal O}}

\usepackage{times}
\usepackage{graphics}
\usepackage{amssymb}
\usepackage{amscd}
\usepackage{amsmath}
\input xy 
\xyoption{all}
\title{Lusztig's partition and sheets}
\author{Giovanna Carnovale}
\newtheorem{theorem}{Theorem}[section]
\newtheorem{lemma}[theorem]{Lemma}
\newtheorem{corollary}[theorem]{Corollary}
\newtheorem{proposition}[theorem]{Proposition}

\newtheorem{remark}[theorem]{Remark}

\newtheorem{lm}[theorem]{Lemma}

\author{Giovanna Carnovale}

\date{With an Appendix by M. Bulois}

\begin{document}
\maketitle
\begin{abstract}
We show that, for a connected reductive algebraic group $G$ over an algebraically closed field of zero or good characteristic, the parts, called strata,  in the partition of $G$ recently introduced by Lusztig are unions of sheets of conjugacy classes. For $G$ simple and adjoint we refine the parametrization of such sheets obtained in previous work with F. Esposito.  We give a simple combinatorial description of strata containing spherical conjugacy classes, showing that  Lusztig's correspondence induces a bijection between unions of spherical conjugacy classes and unions of classes of involutions in the Weyl group. Using ideas from the Appendix by M. Bulois, we show that the closure of a stratum is not necessarily a union of strata.
\end{abstract}



\section{Introduction}

The sheets for the action of an algebraic group on a variety $X$ are the maximal irreducible subsets of $X$ consisting of points whose orbit has fixed dimension. Many important invariants of such actions are preserved along sheets. Sheets for the action of a complex connected reductive algebraic group on its Lie algebra are very well understood \cite{BK, bo}.
Along similar lines, a parametrization and a description of sheets of conjugacy classes in a connected reductive algebraic group $G$ over an algebraically closed field of zero or good characteristic has been given in \cite{gio-espo}.

G. Lusztig  defined in \cite{lustrata} a partition of a connected reductive algebraic group $G$ over an algebraically closed field into certain unions of conjugacy classes of the same dimension. The parts of this partition are called strata. They are the fibers through a map $\phi_G$ from $G$ to a subset of the set ${\rm Irr}(W)$ of isomorphism classes of irreducible representations of the Weyl group $W$ of $G$. On unipotent classes, the map $\phi_G$ coincides with Springer correspondence (with trivial representation of the component group). Lusztig observes that for $G=GL_n(k)$ sheets  coincide with  strata but, for other groups, sheets no longer form a partition and strata are in general not connected. The first result of this paper is Theorem \ref{thm:strata}, stating that in zero or good characteristic every stratum is a union of sheets. In other words, the extension $\phi_G$ of Springer's correspondence is constant along sheets.  This is a direct consequence of the results in \cite{gio-espo} together with compatibility of induction of unipotent conjugacy classes with truncated induction \cite{lusp}. The image of $\phi_G$ contains more irreducible representations of $W$ than those obtained by the Springer correspondence for a trivial local system: this shows once more that, as opposed to the Lie algebra case, where every sheet contains a unique nilpotent orbit \cite{BK}, not every sheet of conjugacy classes contains a unipotent one. As a consequence of Theorem \ref{thm:strata}, we show that strata are locally closed, answering a question of Lusztig. 

Sections \ref{sec:isolated} and \ref{refinement} are devoted to the refinement of some results in \cite{gio-espo}. There, sheets were parametrized by $G$-conjugacy classes of triples $(M, Z(M)^\circ s, \OO)$, where $M$ is the connected centralizer of a semisimple element $s\in G$; $Z(M)^\circ s$ is a suitable coset in $Z(M)/Z(M)^\circ$; and $\OO$ is a rigid unipotent conjugacy class in $M$. A sheet contains a unipotent conjugacy class (up to a central element) if and only if $M$ is a Levi subgroup (of a parabolic subgroup) of $G$ and if this is the case, it is unique. In order to  provide a suitable replacement for the missing unipotent class in a sheet, we show in Proposition \ref{prop:isolated} that every sheet of conjugacy classes contains so-called isolated conjugacy classes (cf. \cite[Definition 2.6]{lusztig-inventiones}). These are finitely many for every semisimple group,  they include unipotent classes and coincide with them (only) if all simple factors in $G$ are of type $A_n$. Isolated classes play a role both in the generalized Springer correspondence and in the representation theory of quantum groups at the roots of unity \cite{DCK}, where they are called exceptional. Isolated classes in a sheet are far from being unique and if two sheets intersect non-trivially, then the intersection contains at least an isolated class. Using injectivity of Springer correspondence we show that two sheets containing a unipotent class meet if and only if they contain the same unipotent class (Proposition \ref{prop:intersection}). 
 
The proof of Theorem \ref{thm:strata} shows that $\phi_G$ depends only on the $G$-class of the pair $(M,\,\OO)$, where $M$ and $\OO$ are as above. This fact suggests that the second term in the triple parametrizing sheets could be dropped. Indeed, Theorem \ref{thm:refinement} states that 
for $G$  simple of adjoint type sheets are parametrized by $G$-classes of such pairs. Since $G$-classes of connected centralizers of semisimple elements are classified in \cite{sommers} and rigid unipotent classes are listed in \cite{spalt}, we obtain a simpler parametrization of sheets in $G$. 

A natural question is  which sheets lie in a given stratum. There is one family of sheets for which the answer is particularly clear:
these are the sheets containing spherical conjugacy classes. We recall that a conjugacy class is called spherical if a(ny) Borel subgroup $B$ of $G$ acts on it with a dense orbit. The property of being spherical is preserved along sheets \cite{arzha}. We show in Section \ref{spherical} that a similar property holds for strata and we describe strata consisting of spherical sheets in combinatorial terms. More precisely, such strata are in bijection with conjugacy classes in the Weyl group $W$ containing a maximum $w_m$, and a spherical conjugacy class $\gamma$ lies in such a stratum if and only if $B w_m B\cap \gamma$ is dense in $\gamma$. This result is a consequence of the combinatorial description of spherical conjugacy classes \cite{ccc, gio-mathZ, cc,lu} and the alternative description of strata in terms of the Bruhat decomposition of $G$ in \cite{lustrata}. Through this alternative description it is proved in Theorem \ref{thm:spherical} that spherical strata correspond to unions of classes of involutions in $W$ having $w_m$ as a maximum. 

In the Appendix by M. Bulois it is shown that, for sheets of adjoint orbits in a Lie algebra, the closure of a sheet is not necessarily a union of sheets. Making use of his counterexamples we show  that, even in the spherical case, the closure of a stratum is not necessarily a union of strata. We give two counterexamples: one for each construction of the strata. 

\subsection{Notation}

Unless otherwise stated, $G$ is a connected, reductive algebraic group over an algebraically closed field $k$ of zero or good characteristic.
Let $T$ be a fixed maximal torus of $G$,  and let $\Phi$ be the associated
root system. Let $B\supset T$ be a Borel subgroup with unipotent radical $U$, let
$\Delta=\{\alpha_1,\ldots,\alpha_n\}$ be the basis of $\Phi$ relative
to $(T,\,B)$.   If $\Phi$ is irreducible, we denote by $-\alpha_0$ the highest positive root in $\Phi$. The Weyl group of $G$ is denoted by $W$,  $\ell$ is the length function on $W$ and ${\rm rk}$ is the rank in the geometric representation of $W$. For $\Pi\subset \Delta$, we denote by $W_\Pi$ the parabolic subgroup of $W$ generated by the simple reflections with respect to roots in $\Pi$, by $w_\Pi$ its longest element, and by $\Phi_\Pi$ be the root subsystem of $\Phi$ generated by $\Pi$. The groups $Sp_{2n}(k)$ and $SO_n(k)$ will fix a bilinear form whose associated matrix with respect to the canonical basis is anti-diagonal.
The numbering of the simple roots is chosen as in \cite{bourbaki}.

Let $G$ act regularly on an irreducible variety $X$. A sheet for this action is an irreducible component of
any of the locally closed subsets $X_{(n)}=\{x\in X~|~ \dim G\cdot x=n\}$, and it is a union of $G$-orbits. For a subset $Y\subset X$, if $m$ is the maximum integer $n$ for which $Y\cap X_{(n)}\neq\emptyset$,  the intersection $X_{(m)}\cap Y$ will be denoted by $Y^{reg}$. Let $V$ be a variety and let $x\in V$. We shall denote by $V_x$ the connected component of $V$ containing $x$ so that, if $V$ is an algebraic group,  we have $V_1=V^\circ$.  When we write $g=su$ for $g\in G$, we mean that $su$ is the Jordan decomposition of $g$, with $s$ semisimple and $u$ unipotent.  If $s\in H^{reg}$ for a subgroup $H$ of $G$, then $(H^{reg})_s$ is well defined, we have $(H^{reg})_s=(H_s)^{reg}$ and we denote it by $H^{reg}_s$. 
 The action of $g\in G$ on an element $x\in G$ by conjugation is indicated by $g\cdot x$. The centralizer of $x$ in $G$ is denoted by $G^x$. Let, for $H$ a connected reductive algebraic group  $\rho^H_y$ be Springer's representation of the Weyl group of $H$ associated with the unipotent element $y\in H$ and trivial local system. If $s\in G$ is semisimple, $W_s$ denotes the Weyl group of $G^{s\circ}$ embedded into $W$ as in \cite{lustrata}. Finally, for finite subgroups $W_1\leq W_2$ of $W$, we denote by ${\bf j}_{W_1}^W$ the $j$-induction functor in \cite[\S 3.2]{lusp}, whenever it is well-defined.

\subsection{Acknowledgements}
G.C. was partially supported by Grants CPDA105885 and CPDA125818/12 of the University of Padova. She thanks Francesco Esposito for pleasant and interesting discussions and the referee for helpful comments and remarks.


\section{Lusztig's strata are union of sheets}

In this section we will show that the parts in Lusztig's partition of $G$ in  \cite{lustrata} are union of sheets. 

We recall that the {\em Jordan class}  $J(g)$ of an element $g=su$ in $G$ is the set $G\cdot((Z(G^{s\circ})^\circ s)^{reg}u)$. Jordan classes were introduced in \cite[\S 3.1]{lusztig-inventiones}, where it is shown that they form a partition of $G$ into locally closed irreducible smooth $G$-stable subsets. In the same paper  the group $L=C_G(Z(G^{s\circ})^\circ)$ for a semisimple element $s\in G$ is introduced. It is the minimal Levi subgroup of a parabolic subgroup containing $G^{s\circ}$. These objects are crucial in the description of sheets.

\begin{theorem}\label{thm:strata} Let $G$ be a connected reductive group in good or zero characteristic. Then, every Lusztig's stratum is a union of sheets.
\end{theorem}
\pf By \cite{gio-espo}, a sheet $S$ of $G$ contains a unique dense Jordan class $J=J(su)$, and, for $L=C_G(Z(G^{s\circ})^\circ)$ we have
$$S=\bigcup_{z\in Z(G^{s\circ})^\circ} G\cdot (sz{\rm Ind}_{L^{zs\circ}}^{G^{zs\circ}}(L^{zs\circ}\cdot u))$$
where ${\rm Ind}$ denotes induction of conjugacy classes as in \cite[\S 1.2]{lusp}. We recall that if $p$ is good for $G$ then it is good for any connected centralizer of a semisimple element in $G$ \cite[Proposition 16]{mcninch-sommers}, and that the algorithm in \cite[II.7]{spalt} for describing induction in classical groups  and the tables for exceptional groups are uniform in good characteristic \cite[p. 176]{spalt}.
On the one hand, $G^{s\circ}\subset C_G(Z(G^{s\circ})^\circ)\cap G^{zs\circ}\subset L^{zs}$ for any $z\in Z(G^{s\circ})^\circ$.
On the other hand, if $x\in L^{sz\circ}\subset L\cap G^{zs\circ}$, then $x$ commutes with $s$ and $zs$, hence $L^{zs\circ}\subset G^{s\circ}$.  Therefore 
$$S=\bigcup_{z\in Z(G^{s\circ})^\circ} G\cdot (sz{\rm Ind}_{G^{s\circ}}^{G^{sz\circ}}(G^{s\circ}\cdot u)).$$
The parts in the partition in \cite{lustrata} are given by the fibers through a map $\phi_G\colon G\to  {\rm Irr}(W)$. 
This is defined on $g=su\in G$ as $\phi_G(g)={\bf j}_{W_s}^W\rho_u^{W_s}$.
We shall compute the image of $x$  lying in a sheet $S$ such that $\overline{S}=\overline{J(su)}$. Since the map $\phi_G$ is constant on conjugacy classes we may assume $x=zsv$ for $v\in {\rm Ind}_{G^{s\circ}}^{G^{sz\circ}}(G^{s\circ}\cdot u)$ and $z\in Z(G^{s\circ})^\circ$.  
Then, $\phi_G(x)={\bf j}_{W_{zs}}^W\rho^{G^{zs\circ}}_v$. By \cite[Theorem 3.5]{lusp}, \cite[\S 6]{lus-indag} we have 
$$\phi_G(x)={\bf j}_{W_{zs}}^W\rho_v^{G^{zs\circ}}={\bf j}_{W_{zs}}^W{\bf j}_{W_{s}}^{W_{zs}}\rho_u^{G^{s\circ}}={\bf j}_{W_s}^W\rho_u^{G^{s\circ}}$$ which depends only on $s$ and $u$, yielding the claim.\hfill$\Box$


\begin{corollary}\label{cor:loc-closed}Lusztig's strata and sheets are locally closed.
\end{corollary}
\pf Let $X$ be a stratum. By Theorem \ref{thm:strata} 
$X=\bigcup_{j=1}^lS_j$ for some sheets $S_j=\overline{J(g_j)}^{reg}$. As $X\subset G_{(n)}$ for some $n$, we have $S_j=\overline{J(g_j)}\cap G_{(n)}$ for every $j$. We recall that
$\overline{J(g_j)}\subset \bigcup_{m\leq n}G_{(m)}=\overline{G_{(n)}}$ so $G_{(>n-1)}:=\bigcup_{m>n-1}G_{(m)}$ is open. 
Then, $X=\bigl(\bigcup_{j=1}^l\overline{J(g_j)}\bigr)\cap G_{(>n-1)}$ is locally closed.
The same argument for $l=1$ proves the result for sheets.
\hfill$\Box$

\section{Isolated elements}\label{sec:isolated}

By \cite[Theorem 5.6]{gio-espo}  the map $S= \overline{J(su)}^{reg}\mapsto (G^{s\circ}, Z(G)Z(G^{s\circ})^\circ s, G^{s\circ}\cdot u)$ induces a bijection between the set of sheets in $G$ and $G$-conjugacy classes of triples $(M,\,Z^\circ Z(G) s,\OO)$ where: $M$ is a {\em pseudo-Levi subgroup of $G$}, i.e., the connected centralizer of a semisimple element $s$ in $G$;   $Z$ is the centre of $M$; the coset $Z^\circ Z(G)s$ s a generator of the cyclic group $Z/Z^\circ Z(G)$; and $\OO$ is a rigid unipotent conjugacy class in $M$. In contrast to the Lie algebra case, where sheets always contain a unique nilpotent orbit, sheets of conjugacy classes do not always contain a unipotent one. Indeed, a sheet contains a unipotent class up to a central element if and only if the term $M$ in the corresponding triple  is the Levi subgroup of a parabolic subgroup. If this is the case, such a class is unique. Following \cite{lusztig-inventiones} we will say that an element $g=su\in G$ is {\em isolated} if $C_G(Z(G^{s\circ})^\circ)=G$, or, equivalently, if $Z(G^{s\circ})/Z(G)$ is finite. Unipotent classes are clearly isolated and, for $G$ simple of type $A_n$, the two definitions coincide.


\begin{proposition}\label{prop:isolated}Every sheet $S= \overline{J(su)}^{reg}$ contains an isolated conjugacy class.
\end{proposition}
%
%
\proof It is enough to prove the statement for $G$ simple of adjoint type. We assume $s\in T$ and that the root system of $G^{s\circ}$ relative to $T$ has a basis $J$ in $\Delta\cup\{\alpha_0\}$. If $J\subset \Delta$ then $G^{s\circ}$ is a Levi subgroup and there is a unipotent conjugacy class in $S$. 
Suppose this is not the case. If $|J|$ equals $r$, the semisimple rank of $G$, then $s$ is isolated and there is nothing to prove.
Otherwise, we consider any $\tilde{J}$ such that $J\subset \tilde{J}\subset \Delta\cup\{\alpha_0\}$ and $|\tilde{J}|=r$. Then, $\tilde{J}$ generates the root system of a pseudo-Levi subgroup $M$ containing $G^{s\circ}$. By construction, $G^{s\circ}$ is a Levi subgroup of $M$ so $Z(G^{s\circ})^\circ s=tZ(G^{s\circ})^\circ$ for some $t\in Z(M)$ \cite[Lemma 33]{mcninch-sommers}. Thus, $G^{t\circ}\supset M$, so  $G=C_G(Z(M)^\circ)\subset C_G(G^{t\circ})$, 
 %
 hence $t$ is an isolated semisimple element in $Z(G^{s\circ})^\circ s$. Thus, any element in $t{\rm Ind}_{G^{s\circ}}^{G^{t\circ}}(G^{s\circ}\cdot u)$ is an isolated element in $S$. \hfill$\Box$
 
\begin{remark}The isolated element in a sheet $S$ is not unique, even up to a central element. For instance, we may consider $G=Sp_{10}(k)$ and the diagonal matrix $s={\rm diag}(-1,a,b,b,1,1,b^{-1},b^{-1},a^{-1},-1)$ with $a\neq b\in k^*\setminus\{\pm1\}$. Then $M:=G^{s\circ}\simeq Sp_2(k)\times SL_2(k)\times Sp_2(k)$ is of type $C_1\times \tilde{A}_1\times {C}_1$ and corresponds to the roots $2(\alpha_1+\alpha_2+\alpha_3+\alpha_4)+\alpha_5$, $\alpha_3$ and $\alpha_5$. Let $S=\overline{J(s)}^{reg}$ be the sheet associated with $(M, s Z(M)^\circ, 1)$. There are $g=rv$ and $h=r'v'$ in $S$ with semisimple parts $r={\rm diag}(-1,-1,1,1,1,1,1,1,-1,-1)$ and $r'={\rm diag}(-1,1,1,1,1,1,1,1,1,1,-1)$. Both $g$ and $h$ are are isolated. More precisely, $G^{r\circ}\simeq Sp_4(k)\times Sp_6(k)$ corresponds to the roots $2(\alpha_1+\alpha_2+\alpha_3+\alpha_4)+\alpha_5$, $\alpha_1$,  $\alpha_3$, $\alpha_4$ and $\alpha_5$ and $G^{r'\circ}\simeq Sp_2(k)\times   Sp_8(k)$ corresponds to the roots $2(\alpha_1+\alpha_2+\alpha_3+\alpha_4)+\alpha_5$, $\alpha_2$, $\alpha_3$, $\alpha_4$ and $\alpha_5$
\end{remark}

\begin{remark}\label{rk:jordan}By \cite[Proposition 48]{gio-espo} the regular part of the closure of any Jordan class $J(su)=G\cdot (Z(G^{s\circ})^\circ s)^{reg}u$ equals
\begin{equation}\label{eq:reg-clo-jo}\overline{J(su)}^{reg}=\bigcup_{z\in Z(G^{s\circ})^\circ}G\cdot zs{\rm Ind}_{G^{s\circ}}^{G^{zs\circ}}(G^{s\circ}\cdot u).\end{equation}  The argument of Proposition \ref{prop:isolated} shows that $\overline{J(su)}^{reg}$ contains isolated elements.
\end{remark}

Any sheet is an irreducible component of the stratum containing it, thus if two sheets have non-empty intersection the stratum containing them is not smooth. It is not hard to see that two sheets in a Lie algebra meet if and only if they contain the same nilpotent orbit. The following proposition is an analogue of this fact. 


\begin{proposition}\label{prop:intersection}If the intersection of two sheets $S_1$ and $S_2$ in $G$ is non-empty, then it contains an isolated class. 
If, in addition, $\OO_1\subset S_1$ and $\OO_2\subset S_2$ for some unipotent classes $\OO_1$ and $\OO_2$, then $\OO_1=\OO_2$.
 \end{proposition}
\pf Let $S_1=\overline{J(g_1)}^{reg}$ and $S_2=\overline{J(g_2)}^{reg}$ be two sheets in $G_{(n)}$ having non-empty intersection and let $G\cdot su\subset S_1\cap S_2$. Since the closure of a Jordan class is a union of Jordan classes (\cite{lusztig-inventiones,gio-espo}), the Jordan class $J(su)$ containing $G\cdot su$ satisfies
$$J(su)\subset S_1\cap S_2=\overline{J(g_1)}\cap \overline{J(g_2)}\cap G_{(n)}$$ and therefore 
$$\overline{J(su)}^{reg}=\overline{J(su)}\cap G_{(n)}\subset S_1\cap S_2=\overline{J(g_1)}\cap \overline{J(g_2)}\cap G_{(n)},$$ so the first statement follows from Remark \ref{rk:jordan}. 
For the second one we observe that if $S_1\cap S_2\neq\emptyset$ then $\phi_G(x)=\phi_G(y)$ for every $x\in\OO_1$, $y\in \OO_2$ and we invoke injectivity of the Springer correspondence on unipotent classes.
\hfill$\Box$

\section{A refinement of the parametrization of sheets}\label{refinement}

The proof of Theorem \ref{thm:strata} shows that the image of $\phi_G$ depends only on the terms $M$ and $\OO$ in the triple corresponding to a sheet. This suggests that the parametrization in \cite{gio-espo} may be improved, and this is in fact the case. We show that the second term in the triple parametrizing sheets  may be dropped when $G$ is simple and of adjoint type. The conjugacy classes of pseudo-Levi subgroups can be deduced from \cite[\S 2.2]{sommers} and rigid unipotent classes are classified in \cite[\S II.7\&II.10]{spalt}, thus a classification of sheets in $G$ follows from these data.  

\begin{theorem}\label{thm:refinement}Let $G$ be simple and of adjoint type. The sheets in $G$ are in bijection with the $G$-conjugacy classes of pairs $(M,\OO)$ where $M$ is a pseudo-Levi subgroup and $\OO$ is a rigid unipotent conjugacy class in $M$.
\end{theorem}
\pf  We show that the $G$-conjugacy class of a triple $(M,\,Z^\circ s,\,\OO)$ where: $M$ is a pseudo-Levi subgroup of $G$ with centre $Z$; the coset $Z^\circ s$ is a generator of the cyclic group $Z/Z^\circ$; and $\OO$ is a rigid unipotent conjugacy class in $M$, is completely determined by the pair $(M,\OO)$.  We may always assume that $s\in T$, so $Z^\circ s\subset T$, and that $M$ is generated by $T$ and by the root subgroups ranging in a subset $\Pi$ of the extended Dynkin diagram, and their opposites. The map $(M, Z^\circ s,\OO)\mapsto (M,\,\OO)$ induces a well-defined and surjective map on the set of $G$-conjugacy classes of triples as above. 

We shall assume $G$ to be of exceptional type because by \cite[Lemma 33]{mcninch-sommers}, $m=|Z/Z^\circ|\leq 2$ in classical groups.

By \cite[Proposition 7]{sommers} any pair of cosets generating $Z/Z^\circ$ are conjugate by some  $w\in W$ preserving $\Pi$, whence preserving $M$. The statement is given when the ground field has characteristic $0$ but the proof holds in good characteristic. 
%

We consider two distinct representatives of elements in the fiber of the $G$-class of $(M,\OO)$. It is not restrictive to assume that they are $(M,Z^\circ s,\OO)$ and $(M,Z^\circ r,\OO)$ with $s\in T$, $r=ws$, and $w$ as above.  
Since $G$ is of exceptional type, $w$ necessarily preserves the unique, if existing, component of type different from type $A$.  Rigid unipotent conjugacy classes in simple algebraic groups are characteristic (\cite[4.5]{bo}), and they are trivial in type $A$. Thus, $w\OO=\OO$ and the induced map is injective. \hfill$\Box$


\section{Spherical sheets and involutions in the Weyl group}\label{spherical}

In this section we shall assume that $G$ is simple for convenience. Theorem~\ref{thm:strata} raises the problem of describing which sheets lie in a given stratum. Using an alternative description of the partition, to be found in \cite[\S5]{lustrata}, we provide a combinatorial answer for strata containing a spherical conjugacy class, showing that they correspond to certain unions of conjugacy classes of involutions in $W$.
We recall that a transitive $G$-space is called {\em spherical} if it has a dense $B$-orbit. We shall denote by ${G}_{sph}$ the union of all spherical conjugacy classes in $G$.
For a conjugacy class $\gamma$ in $G$, let $m_\gamma$ be the uniquely determined Weyl group element such that $\gamma\cap B m_\gamma B$ is dense in $\gamma$. We set $C^\gamma=W\cdot m_\gamma$. By construction, $m_\gamma$ is a maximal length element in $C^\gamma$. By \cite[Theorem 2.13]{cc}, it is also a maximum in $C^\gamma$ with respect to the Bruhat ordering. We have $\dim\gamma\geq \ell(m_\gamma)+{\rm rk}(1-m_\gamma)$ and equality implies that $\gamma$ is spherical, \cite[Theorem 5]{ccc}. 
By construction, if $\gamma\cap BwB\neq\emptyset$ then $w\leq m_\gamma$.
Similarly, for  $S$  a sheet of conjugacy classes, there is a unique element $m_S$ in $W$ such that $B m_S B\cap S$ is dense in $S$. Then, for every $\sigma\in W$ with $B\sigma B\cap S\neq\emptyset$ we 
have $\sigma\leq m_S$ in the Bruhat ordering. Therefore, if $\gamma$ lies in $S$, then $m_\gamma\leq m_S$. It follows from \cite[Lemma 3.1]{clt} that $m^2_{\gamma}=1$ for every class $\gamma$. The same argument shows that $m^2_S=1$ for every sheet $S$. 

It has been shown in \cite{arzha} that, for ${\rm char}(k)=0$, the property of being spherical is preserved along sheets.  As the classification of spherical conjugacy classes in good and odd characteristic \cite{gio-pacific} has the same combinatorics as for ${\rm char}(k)=0$, it follows from the combinatorial description of sheets that the same property holds for conjugacy classes in good and odd characteristic. We will deal now with the case  ${\rm char}(k)=2$ for $\Phi$ of type $A_n$. The result below has already been proved, with different methods, in unpublished work by Mauro Costantini.

\begin{lemma}\label{lem:char2}Assume ${\rm char}(k)=2$. Then the spherical elements in $G=SL_n(k)$ are either involutions up to a scalar or semisimple matrices with at most two eigenvalues. For a spherical conjugacy class $\gamma$ we have $\dim \gamma=\ell(m_\gamma)+{\rm rk}(1-m_\gamma)$ and the property of being spherical is constant along sheets.
\end{lemma}
\pf The argument in \cite[Proposition 1]{arzha} shows that if for a sheet $S=\overline{J(g)}^{reg}$ we have $S\cap G_{sph}\neq\emptyset$ then $J(g)\subset G_{sph}$.
In $SL_n(k)$ the Jordan classes that are dense in a sheet are precisely those consisting of semisimple classes. The only spherical semisimple elements in $SL_n(k)$ are those with at most two eigenvalues \cite[Table 1]{kr}. Therefore the only non-semisimple spherical classes lie in $\overline{J(s)}^{reg}$ for some semisimple element $s$ with at most two distinct eigenvalues. Such sheets contain only semisimple elements and unipotent elements (up to a scalar) corresponding to a partition of type $(2^a,1^{n-2a})$. The latter are spherical and the dimension formula holds for them \cite[3.1.1]{mauro-cattiva}. The dimension formula for spherical semisimple classes follows from a direct computation as in \cite[Theorem 15]{ccc}. 
\hfill$\Box$

%
%
\bigskip
By abuse of notation the sheets contained in $G_{sph}$ will be called {\em spherical sheets}. We will prove that $m_\gamma$ is constant along spherical sheets. 
 
\begin{lemma}\label{involutions}Let $w,\,\sigma$ be two involutions in $W$ such that $w\leq \sigma$ and $\ell(w)+{\rm rk}(1-w)=\ell(\sigma)+{\rm rk}(1-\sigma)$. Then $w=\sigma$.
\end{lemma}
\pf By \cite{incitti1,incitti2,incitti3} for the classical groups and  \cite[Theorems 4.2, 4.8]{hult1} in the general case, the poset of involutions in a Weyl group is graded with rank function $\rho(w)=\frac{1}{2}(\ell(w)+\ell^a(w))$ where the absolute length $\ell^a$ of $w$ is the minimal number of reflections in $W$ needed to  express $w$ as a product of reflections. Thus, if $w\leq\sigma$ then $\rho(w)\leq\rho(\sigma)$ and equality holds only if $w=\sigma$.
By a result of Kostant in \cite{kostant}  we have $\ell^a(w)={\rm rk}(1-w)$, whence the statement.  \hfill$\Box$


\begin{proposition}\label{pro:constant}Let $S$ be a spherical sheet. Then, for every conjugacy class $\gamma$ lying in $S$ we have $m_\gamma=m_S$.
\end{proposition}
\pf Let $\gamma, \gamma'$ be  conjugacy classes in $S$, with $\gamma'\cap Bm_S B\neq\emptyset$. Then, $m_{\gamma'}=m_S$ because $m_{\gamma'}$ is maximal among the Weyl group elements whose Bruhat double coset meets $\gamma'$. Since $\gamma$ and $\gamma'$ are spherical
we have $\dim(\gamma)=\ell(m_\gamma)+{\rm rk}(1-m_\gamma)$ and $\dim(\gamma')=\ell(m_{\gamma'})+{\rm rk}(1-m_{\gamma'})$  by \cite{ccc, lu, gio-mathZ} if ${\rm char}(k)\neq2$ and Lemma \ref{lem:char2} if  ${\rm char}(k)=2$. Since $\gamma$ and $\gamma'$ lie in the same sheet we have
 $$\ell(m_\gamma)+{\rm rk}(1-m_\gamma)=\ell(m_{\gamma'})+{\rm rk}(1-m_{\gamma'})$$ and $m_\gamma\leq m_{\gamma'}=m_S$. 
%
%
Lemma \ref{involutions} applies.
\hfill$\Box$

\begin{remark}For $k={\mathbb C}$ and $\gamma$ a spherical conjugacy class,  $m_\gamma$ is strictly related to the $G$-module decomposition of ${\mathbb C}[\gamma]$. Indeed,  it is well-known that ${\mathbb C}[\gamma]$  is multiplicity-free. In addition, the highest weights occurring with multiplicity $1$ generate a finite index sub-lattice among those integral weights $\lambda$ such that $m_\gamma \lambda=-\lambda$ and $-w_0\lambda=\lambda$  (\cite{ccc,mauro-mathZ}). Broadly speaking, Proposition~\ref{pro:constant} may be seen as a discrete analogue to \cite[Theorems 3.5, 3.8]{BK} for spherical conjugacy classes.
\end{remark}

We recall the alternative approach to strata in \cite[\S2]{lustrata}. The $G$-orbits of pairs of Borel subgroups in $G$ are parametrized by the elements of $W$. We denote such orbits by $\OO_w$. 
For $w\in W$, let $$G_w=\{g\in G~|~(B', gB'g^{-1})\in \OO_w,\,\mbox{ for some Borel subgroup }B' \mbox{ of $G$}\}.$$
In other words, $G_w$ is the union of all conjugacy classes $\gamma$ in $G$ such that $\gamma\cap BwB\neq\emptyset$.
For $C$ a conjugacy class in $W$, let $C_{min}$ ($C_{max}$, respectively) denote its set of minimal length elements (maximal length elements, respectively). 
For $w,w'\in C_{min}$ we have $G_w=G_{w'}$ by \cite[1.2(a)]{lusztig1}  and \cite[8.2.6(b)]{GP}. We denote by $G_C$ the set $G_w$ for $w\in C_{min}$. Let $\delta_C$ be the minimal dimension of a conjugacy class $\gamma$ contained in $G_C$ and let $\underline{G}_C$ be the union of all classes in $G_C$ of dimension exactly $\delta_C$. According to \cite[Theorem 5.2]{lustrata}, whose proof is announced for classical groups and explicit for exceptional groups, the set $\underline{G}_C$ is a stratum and all strata can be described this way. 
 
Let ${W}_{inv}$ be the set of involutions in $W$ and for a conjugacy class $C$ in $W$, let $\underline{W}_C$ be the union of conjugacy classes $C'$ in $W$ such that $\underline{G}_C=\underline{G}_{C'}$. We set
$${\mathcal T}:=\{\Pi\subset \Delta~|~w_0(\alpha)=w_\Pi(\alpha),\,\forall\alpha\in\Phi_\Pi\}.$$
For $C$ a class in ${W}_{inv}$, all elements in $C_{max}$ are of the form $ww_\Pi$ for some $\Pi\in{\mathcal T}$, \cite[Theorem 1.1(ii)]{per-row}. For $\Pi,\Pi'\in {\mathcal T}$ we have $w_0w_\Pi\leq w_0 w_{\Pi'}$  if and only if $\Pi\supset\Pi'$. We also set
$${\mathcal M}:=\{\Pi\in{\mathcal T}~|~ w_0w_\Pi\mbox{ is the unique maximal length elements in its $W$-class}\}.$$

\begin{lemma}\label{lem:subsets}Let $\Pi\in{\mathcal T}\setminus{\mathcal M}$. Then 
\begin{enumerate}
\item The set of elements  $\Pi'$ in ${\mathcal M}$  satisfying $\Pi'\subset\Pi$ has a maximum $M_\Pi$ with respect to inclusion.
\item $\Pi$ is the union of $M_\Pi$ and some isolated simple roots orthogonal to $M_{\Pi}$.
\end{enumerate}
\end{lemma}
\pf The list of elements in ${\mathcal M}$ is given in \cite[Lemma 3.5]{clt}. A straightforward verification gives 1. and 2.
\hfill$\Box$

\begin{lemma}\label{lem:intersections}Let $\gamma$ be a conjugacy class in $G$ such that $\gamma\cap B w_0w_\Sigma B\neq\emptyset$ for some $\Sigma\in {\mathcal T}$, and let $\Pi=\Sigma\cup\Sigma'\in\mathT$ for some $\Sigma'= \{\beta_1,\,\ldots,\beta_l\}\subset \Delta$ with $(\beta_i ,\beta)=0$ for every $i$ and for every $\beta\in \Pi\setminus\{\beta_i\}$. Then 
$\gamma\cap B w_0w_\Pi B\neq\emptyset$. \\ In particular, if $\gamma \cap B w_0w_{M_\Pi}  B\neq\emptyset$ for some $\Pi\in{\mathcal T}$, then $\gamma\cap B w_0w_{\Pi} B\neq\emptyset$.
\end{lemma}
\pf The proof is by induction on $l$, the case of $l=0$ being trivial. Assume the statement is proved for $l=i$. Let $\Sigma_i=\Sigma\cup\{\beta_1,\,\ldots,\,\beta_i\}$, $\alpha=\beta_{i+1}$ and assume $\gamma\cap B w_0w_{\Sigma_i}B\neq\emptyset$. Then, there exists $x= \w_0\w_{\Sigma_i} x_{\alpha}(t)\v\in \gamma\cap \w_0\w_{\Sigma_i}U\cap \gamma$ for some $\w_0 \w_{\Sigma_i}\in N(T)$ representing  $w_0 w_{\Sigma_i}$, some $t\in k$ and some $\v\in P_\alpha^u$, the unipotent radical of the minimal parabolic subgroup of $G$ associated with $\alpha$. Assume that the parametrization of the root subgroup $x_\alpha(k)$ is chosen as in \cite[Lemma 8.1.4]{springer}. 
There is $\eta\in k^*$ such that  $x_{-\alpha}(\xi)\w_0 \w_{\Sigma_i}=\w_0 \w_{\Sigma_i}x_{\alpha}(\eta \xi)$ for every $\xi\in k$. We choose $\xi\in k$ satisfying  $\eta \xi^2+t\xi-1=0$.
Then, for $y:=x_{-\alpha}(\xi)x x_{-\alpha}(-\xi)\in \gamma$ and $\v'=x_{-\alpha}(\xi)\v x_{-\alpha}(-\xi)\in P_\alpha^u$ we have
$$\begin{array}{ll}
y&=\w_0 \w_{\Sigma_i}x_{\alpha}(\eta \xi+t)\v x_{-\alpha}(-\xi)=\w_0 \w_{\Sigma_i}x_{\alpha}(\xi^{-1})x_{-\alpha}(-\xi)\v'\\
&\in w_0 w_{\Sigma_i} s_\alpha T x_{\alpha}(-\xi^{-1})\v'\subset B w_0w_{\Sigma_{i+1}}B
\end{array}$$
where we have used \cite[Lemma 8.1.4(22)]{springer}. Last statement follows from Lemma \ref{lem:subsets} (2).
\hfill$\Box$

\begin{lemma}\label{lem:w^2}Let $C$ be a class in $W$ and $\gamma$ be a spherical class in $SL_n(k)$ such that $\gamma\subset\underline{G}_C$. Then 
\begin{enumerate}
\item $C\subset{W}_{inv}$
\item $B w B\cap \gamma\neq\emptyset$ for every $w\in C$.
\end{enumerate}
\end{lemma}
\pf 1. For ${\rm char}(k)\neq2$ this is \cite[Theorem 2.7]{gio-mathZ}. If ${\rm char}(k)=2$, then $\Phi$ is of type $A$ and spherical classes are described in Lemma \ref{lem:char2}. If $\gamma$ is the class of an involution there is nothing to prove.  Let thus $\gamma$ be a semisimple class in $SL_n(k)$ with two eigenvalues of multiplicity $m$ and $q=n-m$, respectively, for $m\geq q$. Let $w\in C_{min}$, so $BwB\cap\gamma\neq\emptyset$. If $w$ has no fixed points (elliptic case), we may take $w=(1,2,\ldots,i_1)(i_1+1,\ldots,i_1+i_2)\cdots(i_1+i_2+\cdots+i_{r-1},\ldots,n)$. Then \cite[Lemma 4.1]{clt} gives $r\geq m\geq\left[\frac{n}{2}\right]$, forcing $i_j\leq 2$ for every $j$. 

Assume now that the set of fixed points of $w$ is $K=\{k_1,k_1+k_2,\ldots,k_1+\cdots+k_t\}$, i.e., 
$w$ lies in the parabolic subgroup of $W$ isomorphic to $S_{k_1-1}\times\cdots \times S_{k_t-1}\times S_{n-k_1-\cdots-k_t}$, where some of the factors are possibly trivial. Arguing as in \cite[1.1]{lusztig1}, see also \cite[Theorem 5.2]{sev} for different notation, we see that if $\gamma$ has minimal dimension in $G_C$ then $\gamma \cap L\cap B_L w B_L\neq\emptyset$, where $L$ is the standard Levi subgroup of a standard parabolic subgroup associated with the simple roots indexed by $\{1,\ldots,n\}\setminus \{k_1,k_1+1,k_1+k_2-1,k_1+k_2,\ldots k_1+\cdots k_t-1, k_1+\cdots + k_t\}$ and $B_L=L\cap B$. Then, 
$L=Z(L)^\circ L_1\cdots L_{t+1}$ where $L_j\simeq SL_{k_j}(k)$ and some of the factors are possibly trivial. We work componentwise. As each component of $w$ has no fixed points, we may reduce to the elliptic case.

2. If ${\rm char}(k)=2$, $\Phi$ of type $A$ and $\gamma$ is semisimple,  $C\subset {W}_{inv}$ by 1, so \cite[Theorem 4.2]{clt} applies.
In all other cases \cite[Lemma 2.2]{cc} applies. \hfill$\Box$

\begin{theorem}\label{thm:spherical}Let $C$ be a conjugacy class in $W$, and let $\underline{W}_C$ and $\underline{G}_C$ be as above.
\begin{enumerate}
\item If $\gamma\subset G_{sph}$, then $\gamma\subset \underline{G}_{C^\gamma}\subset G_{sph}\cap \Bigl(\bigcup_{\gamma'\subset G\atop m_{\gamma'}=m_C}\gamma'\Bigr)$ where $\gamma'$ runs through the conjugacy classes in $G$.
\item If $C$ has a maximum $m_C$, 
then 
\begin{equation}\label{GC}\underline{G}_C=G_{sph}\cap \Bigl(\bigcup_{\gamma\subset G\atop m_\gamma=m_C}\gamma\Bigr)=G_{sph}\cap\Bigl(\bigcup_{S\subset G\atop m_S=m_C}S\Bigr)\end{equation} where the $\gamma$'s are conjugacy classes and the $S$'s are sheets in $G$.
\item If $\underline{G}_C\cap G_{sph}\neq\emptyset$  then $C\subset {W}_{inv}$.
\item If $C\subset{W}_{inv}$ then $\underline{G}_C=\underline{G}_{C^\gamma}\subset{G}_{sph}$, for some class $\gamma$  and $m_\gamma=w_0w_{M_\Pi}$ (notation as in Lemma \ref{lem:subsets}) for one (hence for every) $w_0w_\Pi\in C_{max}$. 
\item If $\underline{G}_C\cap G_{sph}\neq\emptyset$  then ${ \underline{W}_C}$ has a maximum which equals $m_\gamma$ for every $\gamma\subset\underline{G}_C$. 
\item If $C$ is a class with a maximum $m_C$, then
\begin{equation}\label{WC}\underline{W}_C=W_{inv}\cap \Bigl(\bigcup_{{C\subset W\atop m_C=w_0w_{M_\Pi}}\atop 
{\rm for }\;w_0w_\Pi\in C'_{max}}C'\Bigr).\end{equation}
\end{enumerate} 
\end{theorem}
\pf 1. Certainly $\gamma \cap B m_\gamma B\neq\emptyset$. Since $\gamma$ is spherical, $B \sigma B\cap\gamma\neq\emptyset$ for every $\sigma\in C^\gamma$, \cite[Lemma 2.2]{cc},  \cite[Theorem 4.2]{clt}. Thus $\gamma\subset G_{C^\gamma}$.  Let $\gamma'\subset \underline{G}_{C^\gamma}$. By \cite[Propositions 2.8, 2.9]{clt},  (a reformulation of \cite[\S 2.9]{GKP} ,\cite[Proposition 5.3.4]{EG}), we have $\gamma'\cap B m_\gamma B\neq\emptyset$. Therefore, $\dim \gamma'\geq \ell(m_\gamma)+{\rm rk}(1-m_\gamma)=\dim \gamma$, where the equality on the right follows from the main result in \cite{ccc,gio-mathZ,lu} and Lemma \ref{lem:char2}.
Hence, $\dim \gamma'=\dim \gamma$,  $\gamma\subset\underline{G}_{C^\gamma}$, $\gamma'$ is spherical by  \cite[Theorem 5]{ccc} and $m_{\gamma'}=m_{\gamma}$. 

2. We claim that for every class $C$ with a maximum $m_C$ there always exists a spherical conjugacy class $\gamma_0$ such that $m_{\gamma_0}=m_C$. 
If ${\rm char} (k)\neq 2$ this is \cite[Remark 3]{clt}. If ${\rm char} (k)=2$ then $\Phi$ is of type $A$. In this case the classes in $W$ having a maximum coincide with the classes of involutions, and the correspondence $\gamma\mapsto m_\gamma$ is a bijection between the set of spherical unipotent classes and the set of classes of involutions in $W$. Hence, the first inclusion $\subset$ follows from 1. On the other hand, if $\gamma'$ is spherical and $m_{\gamma'}=m_C$ then $C=C^{\gamma'}$ and again by 1., we have $\gamma'\subset \underline{G}_{C}$ and the first equality of sets follows. Combining with Proposition \ref{pro:constant} yields the second one.


3. If $\gamma$ is spherical and ${\rm char}(k)\neq2$ then $\gamma\subset \bigcup_{w^2=1}BwB$ by \cite[Theorem 2.7]{gio-mathZ}. For ${\rm char}(k)=2$, $\Phi$ is of type $A$ and we invoke Lemma \ref{lem:w^2}(1).

4. Let $\gamma\subset \underline{G}_{C}$ so $\gamma\cap B w B\neq\emptyset$ for some $w\in C_{min}$. By \cite[Propositions 2.8, 2.9]{clt}, we have $\gamma\cap B \sigma B\neq\emptyset$ for some $\sigma\in C_{max}$. Then $\sigma=w_0w_\Pi$ for some $\Pi\in{\mathcal T}$. 
If $\Pi\in{\mathcal M}$ this is statement 3. so we may assume $\Pi\not\in {\mathcal M}$. Let $M_\Pi$ as in Lemma \ref{lem:subsets}.
We have $\sigma\leq w_0w_{M_\Pi}\leq m_\gamma$ and so
\begin{equation}\label{dim}\ell(w_0w_{M_\Pi})+{\rm rk}(1-w_0w_{M_\Pi})\leq \ell(m_\gamma)+{\rm rk}(1-m_\gamma)\leq \dim \gamma.\end{equation}

Let $C'=W\cdot w_0w_{M_\Pi}$ and let $\gamma'\subset\underline{G}_{C'}$. By Lemma  \ref{lem:intersections} we have $\gamma'\cap Bw_0w_{\Pi}B\neq\emptyset$. On the other hand $\gamma'$ is spherical by statement 2., so $\gamma'\cap B w B\neq\emptyset$, by \cite[Lemma 2.2]{cc} and Lemma \ref{lem:w^2} (2). Therefore $\gamma'\subset G_C$ and 
\begin{equation}\label{dim'}\ell(w_0w_{M_\Pi})+{\rm rk}(1-w_0w_{M_\Pi})=\dim \gamma'\geq\dim \gamma.\end{equation} Thus, the inequalities in \eqref{dim} and \eqref{dim'} are equalities, 
$\gamma\subset G_{sph}$, $m_\gamma=w_0w_{M_\Pi}$, and $\underline{G}_{C^\gamma}\subset \underline{G}_C$. 
By 2., $\underline{G}_C=\underline{G}_{C^\gamma}$, whence the statement. 

5. By 3. and 4, $\underline{G}_C=\underline{G}_{C^\gamma}\subset G_{sph}$, so $m_\gamma\in \underline{W}_C$. Therefore 
it is enough to show that $\gamma\cap BwB\neq\emptyset$ for every $w\in W_C$: this is Lemma \ref{lem:w^2} in type $A$ and \cite[Lemma 2.2]{cc} otherwise. 

6. ($\subset$). If $C'\subset \underline{W}_C$ then by 2. and 5. we have 
$$G_{sph}\cap \bigl(\bigcup_{m_\gamma=m_C}\gamma\bigr)=\underline{G}_{C}=\underline{G}_{C'}=G_{sph}\cap \Bigl(\bigcup_{m_\gamma=w_0w_{M_\Pi}\atop
{\rm for }w_0w_\Pi\in C'_{max}}\gamma\Bigr)$$ 
whence the first inclusion. ($\supset$). If $C'\subset W_{inv}$ and $m_C=w_0w_{M_\Pi}$ for $w_0w_\Pi\in C'_{max}$, then by the argument in 4., for $\gamma\subset\underline{G}_C$ we have $\underline{G}_{C'}=\underline{G}_{C^\gamma}=\underline{G}_{C}$ so $C'\in\underline{W}_C$.\hfill$\Box$

\bigskip 

Let $(G/\!\sim)$ denote the set of strata of the form $\underline{G}_C$ and $(W/\!\sim)$ denote the set of subsets $\underline{W}_C$ of $W$. Theorem \ref{thm:spherical} implies the following fact.

\begin{corollary}Lusztig's bijection induces bijections $$({G}_{sph}/\!\sim) \longleftrightarrow\mathM  \longleftrightarrow({W}_{inv}/\!\sim).$$ 
where the correspondence $\underline{G}_C \leftrightarrow\underline{W}_C$ is given by \eqref{GC} and \eqref{WC}.
\end{corollary}

\begin{remark}The closure of a stratum is in general not a union of strata, not even in the case of spherical strata. We provide 2 counterexamples, stemming from the counterexamples in the Appendix. The first one uses the description of the partition of $G$ in terms of the Bruhat decomposition, the second one uses the description in terms of the map $\phi_G$.
\begin{enumerate}
\item Let $G=SO_{8}(k)$, and let $X$ be the spherical stratum corresponding to $w_0$ as in Theorem \ref{thm:spherical}. By the classification in  \cite{ccc,gio-pacific}, $X$ is the union of $3$ classes: the rigid unipotent class $\OO_1$ with partition $[3,2^2,1]$; $\OO_1$ multiplied by the non-trivial central element $-I$ in $G$; and the conjugacy class of an orthogonal diagonal matrix $s={\rm diag}(1,1,-1,-1,-1,-1,1,1)$. In other words, it is the union of the sheets corresponding to the triples $(G, 1, \OO_1)$, $(G,-I,\OO_1)$ and $(M, s, 1)$, where $M$ is the pseudo-Levi of type $D_2\times D_2$ corresponding to the  simple roots $\alpha_1,\alpha_1+2\alpha_2+\alpha_3+\alpha_4$, $\alpha_3$ and $\alpha_4$. Hence, $\overline{X}\setminus X$ consists only of unipotent classes, up to a central element. In particular, this set contains the unipotent class $\OO_2$ corresponding to the partition $[3,1^5]$, which is spherical and not rigid. Then $\OO_2$ lies in a non-trivial spherical sheet, hence in a non-trivial stratum which cannot be contained in $\overline{X}$. 
\item Let $G=SL_n(k)$. Then sheets coincide with strata by \cite[1.16]{lustrata} and we may use counterexample (2) in the Appendix. 
\end{enumerate}
\end{remark}

\noindent{\sc G. Carnovale}\\
Dipartimento di Matematica - Universit\`a degli Studi di Padova\\
via Trieste 63 - 35121 Padova - Italy\\
email: carnoval@math.unipd.it 

\section*{Appendix by Micha\"el Bulois}\label{appendix}

In this Appendix we answer to a frequently asked question. We focus on the case of sheets for the adjoint action of a semisimple group $G$ on its Lie algebra. We give two families of examples of sheets whose closure is not a union of sheets in this setting.

Let $\mathfrak{g}$ be a semisimple Lie algebra defined over an algebraically closed field $k$ of characteristic zero.
Let $G$ be the adjoint group of ${\mathfrak g}$. For any integer $m$, one defines $${\mathfrak g}_{(m)}=\{x\in\g\mid \dim G\cdot x=m\}.$$
In this case a sheet is an irreducible component of $\g_{(m)}$ for some $m\in{\mathbb N}$. We refer to \cite{BK,bo} for elementary properties of sheets. An important one is that each sheet contains a unique nilpotent orbit.

There exists a well known subdivision of sheets which forms a stratification. The objects considered in this subdivision are Jordan classes and generalize the classical Jordan's block decomposition in ${\mathfrak g}{\mathfrak l}_n$. These classes and their closures are widely studied in \cite{bo}.  Since sheets are locally closed, a natural question is then the following. 
$$\mbox{If $S$ is a sheet, is $\overline{S}$  a union of sheets?}$$
The answer is negative in general. We give two families of counterexamples below.

\begin{enumerate}
\item  A nilpotent orbit $\Od$ of $\mathfrak g$ is said to be rigid if it is a sheet of $\mathfrak g$. Rigid orbits are key objects in the description of sheets given in \cite{bo}. They are classified  in \cite[\S II.7\&II.10]{spalt}.
The closure ordering of nilpotent orbits (or \emph{Hasse diagram}) can be found in \cite[\S II.8\&IV.2]{spalt}. 
One easily checks from these classifications that there may exists some rigid nilpotent orbit $\Od_{1}$ that contains a non-rigid nilpotent orbit $\Od_{2}$ in its closure.
Then, we set $S=\Od_{1}$ and we get $\Od_{2}\subset\overline{S}\subset\NN({\mathfrak g})$ where $\NN({\mathfrak g})$ is the set of nilpotent elements of $\mathfrak g$. 
Since $\Od_{2}$ is not rigid, the sheets containing $\Od_{2}$ are not wholly included in $\NN({\mathfrak g})$. 
Therefore, the closure of $S$ is not a union of sheets.

Here are some examples of such nilpotent orbits.
In the classical  cases, we embed $\mathfrak g$ in ${\mathfrak g}{\mathfrak l}_n$ in the natural way. Then, we can assign to each nilpotent orbit $\Od$, a partition of $n$, denoted by $\Gamma(\Od)$. This partition defines the orbit $\Od$, sometimes up to an element of $\Aut({\mathfrak g})$.
In the case ${\mathfrak g}=\so_{8}$ (type D$_{4}$), there is exactly one rigid orbit $\Od_{1}$, such that $\Gamma(\Od_1)=[3,2^2,1]$. It contains in its closure the non-rigid orbit $\Od_{2}$ such that $\Gamma(\Od_2)=[3,1^5]$. 
Very similar examples can be found in types C and B. 

In the exceptional cases, we denote nilpotent orbits by their Bala-Carter symbol as in \cite{spalt}. 
Let us give some examples of the above described phenomenon. 
\begin{itemize}
\item in type $E_6$ ($\Od_1=3A_1$ and $\Od_{2}=2A_1$), 
\item in type $E_7$ ($\Od_1=A_2+2A_1$ and $\Od_{2}=A_2+A_1$), 
\item in type $E_8$ ($\Od_1=A_2+A_1$ and $\Od_{2}=A_2$), 
\item and in type $F_4$ ($\Od_1=A_2+A_1$ and $\Od_{2}=A_2$).   
\end{itemize}

\item In the case ${\mathfrak g}={\mathfrak s}{\mathfrak l}_{n}$ of type $A$, there is only one rigid nilpotent orbit, the null one. Hence the phenomenon depicted in 1 can not arise in this case. Let $S$ be a sheet and let $\lambda_S=(\lambda_{1}\geqslant\dots\geqslant \lambda_{k(\lambda_S)})$ be the partition of $n$ associated to the nilpotent orbit $\Od_S$ of $S$ according to the size of the blocks in the Jordan form of an element of $\Od_S$. Let $\tilde{\lambda}$  be the dual partition of $\lambda$, i.e. $\tilde{\lambda_{i}}=\#\{j\mid \lambda_{j}\geqslant i\}$ (see, e.g.,  \cite[\S2.2]{Kr})  and let  $\lf_{S}$ be the standard Levi subalgebra whose size of the blocks are the parts of $\tilde{\lambda}_S$.

As a consequence of the theory of induction of orbits, cf. \cite{bo}, we have \begin{equation}\label{eq}\overline{S}=\overline{G\cdot{\mathfrak h}_{S}}^{reg}\end{equation} where ${\mathfrak h}_{S}$ is the centre of $\lf_{S}$.  In particular, the map sending a sheet $S$ to its nilpotent orbit $\Od_S$ is a bijection.

An easy consequence of \eqref{eq} is the following (see \cite[Satz 1.4]{Kr}). Given any two sheets $S$ and $S'$ of $\mathfrak g$, 
we have $S\subset \overline{S'}$ if and only if ${\mathfrak h}_{S}$ is $G$-conjugate to a subspace of ${\mathfrak h}_{S'}$  or, equivalently, $\lf_{S'}$ is conjugate to a subspace of $\lf_{S}$. This can be translated in terms of partitions by defining a partial ordering on the set of partitions of $n$ as follows.
We say that $\lambda\preceq\lambda'$ if there exists a partition $(J_i)_{i\in [\![1,p(\lambda)]\!]}$ of $[\![1,p(\lambda')]\!]$ such that $\tilde{\lambda}_i=\sum_{j\in J_i} \tilde{\lambda}'_j$. Hence, we have the following characterization.
\begin{lm}
$S\subset \overline{S'}$ if and only if $\lambda_S\preceq\lambda_{S'}$.
\end{lm}
One sees that this criterion  is strictly stronger than the characterization of inclusion relations of closures of nilpotent orbits (see, e.g., \cite[\S6.2]{CM}).
More precisely, one easily finds two sheets $S$ and $S'$ such that $\Od_{S}\subset\overline{\Od_{S'}}$ while $\lambda_{S}\npreceq\lambda_{S'}$.
Then, $\Od_S\subset \overline{S'}$, $S$ is the only sheet containing $\Od_S$ and $S\not\subset \overline{S'}$.
For instance, take $\lambda_{S'}=[2,2]$, $\lambda_{S}= [2,1,1]$. Their respective dual partitions being $[2,2]$ and $[3,1]$, we have $\lambda_{S}\npreceq\lambda_{S'}$.
\end{enumerate}

\noindent{\sc Micha\"el Bulois}\\
Universit\'e de Lyon, CNRS UMR 5208,\\
Universit\'e Jean Monnet, Institut Camille Jordan, \\
Maison de l'Universit\'e, 10 rue Tr\'efilerie, CS 82301, \\
42023 Saint-Etienne Cedex 2, France. \\
EMAIL: michael.bulois@univ-st-etienne.fr

%





\begin{thebibliography}{1}

\bibitem{arzha}{\sc I. V. Arzhantsev,}
\newblock{\em On $SL_2$-actions of complexity one,}
\newblock Izv. Ross. Akad. Nauk Ser. Mat. 61(4), 3--18 (1997);
translation in Izv. Math. 61(4) 685--698 (1997).  


\bibitem{bo}{\sc W. Borho}
\newblock{\em \"Uber Schichten halbeinfacher Lie-Algebren,}
\newblock Invent. Math., 65, 283--317 (1981/82). 

\bibitem{BK}{\sc W. Borho, H. Kraft,}
\newblock {\em \"Uber Bahnen und deren Deformationen bei linearen
  Aktionen reduktiver Gruppen,}
\newblock Comment. Math. Helvetici, 54, 61--104  (1979).


\bibitem{bourbaki}{\sc N.\ Bourbaki,}
\newblock{\em \'El\'ements de Math\'ematique. Groupes et Alg\`ebres de Lie, Chapitres 4, 5, et 6,}
\newblock Masson, Paris (1981).



\bibitem{ccc}{\sc N.\ Cantarini, G.\ Carnovale, M.\ Costantini,}
\newblock {\em Spherical orbits and representations of $\Ue$,}
\newblock Transformation\  Groups, 10, No.\ 1, 29--62  (2005).

\bibitem{gio-mathZ}{\sc G.\ Carnovale,}
\newblock {\em Spherical conjugacy classes and involutions in the Weyl
group,} 
\newblock Math. Z. 260(1)  1--23 (2008).

\bibitem{gio-pacific}{\sc G.\ Carnovale,}
\newblock {\em A classification of spherical conjugacy classes in good characteristic,}
\newblock Pacific Journal of Mathematics  245(1) 25--45 (2010).


\bibitem{cc}{\sc  G.\ Carnovale, M.\ Costantini,}
\newblock{\em On Lusztig's map for spherical unipotent conjugacy classes,}
\newblock Bull. London Math. Soc. doi:10.1112/blms/bdt048 . 

\bibitem{gio-espo}{\sc G.\ Carnovale, F.\ Esposito,}
\newblock{\em On sheets of conjugacy classes in good characteristic,}
\newblock IMRN, doi: 10.1093/imrn/rnr047  2012(4) (2012), 810-828.

\bibitem{clt}{\sc K. Y. Chan, J-H. Lu, S. K-M. To,}
\newblock{\em On intersections of conjugacy classes and Bruhat cells,}
\newblock Transform. Groups 15(2), 243--260 (2010).

\bibitem {CM}{\sc D.~H.~Collingwood and W.~M.~McGovern,}
\newblock \emph{Nilpotent orbits in semisimple Lie algebras}, 
\newblock Van Nostrand Reinhold Mathematics Series, New York, (1993).

\bibitem{mauro-mathZ}{\sc M. Costantini,}
\newblock{\em On the coordinate ring of  spherical conjugacy classes,}
\newblock Math. Z. 264  327--359 (2010).

\bibitem{mauro-cattiva}{\sc M. Costantini,}
\newblock{\em A classification of  unipotent spherical conjugacy classes in bad characteristic,}
\newblock  Trans. Amer. Math. Soc. 364(4) 1997Ð2019 (2012).


\bibitem{DCK}{\sc C.\ De Concini, V.\ G.\ Kac,}
\newblock {\em Representations of quantum groups at roots of $1$:
  reduction to the exceptional case,}
\newblock Infinite analysis, Part A, B (Kyoto,1991), 141--149,
Adv. Ser. Math. Phys., 16, World Sci. Publ., River Edge, NJ, (1992).  

\bibitem{EG}{\sc E. Ellers, N. Gordeev, }
\newblock{\em Intersections of conjugacy classes with Bruhat cells in Chevalley groups, }
\newblock Pacific J. Math. 214 no. 2, 245--261, (2004).



\bibitem{GKP}{\sc M. Geck, S. Kim, G. Pfeiffer,}
\newblock{\em Minimal length elements in twisted conjugacy classes
of  finite Coxeter groups, }
\newblock J. Algebra 229  no. 2, 570--600 (2000).

\bibitem{GP}{\sc M. Geck, G. Pfeiffer,}
\newblock{\em Characters of  Finite Coxeter Groups and Iwahori-Hecke Algebras,}
\newblock London Mathematical Society Monographs, New Series 21, Oxford University Press, London,  (2000).


\bibitem{incitti1}{\sc F. Incitti,}
\newblock{\em The Bruhat order on the involutions of the symmetric group,}
\newblock J. Algebraic Combin. 20, 243--261 (2004).

\bibitem{incitti2}{\sc F. Incitti,}
\newblock{\em The Bruhat order on the involutions of the hyperoctahedral group,}
\newblock European J. Combin. 24, 825--848 (2003).

\bibitem{incitti3}{\sc F. Incitti,}
\newblock{\em Bruhat order on the involutions of the  classical Weyl groups,}
\newblock Advances in Applied Mathematics 37, 68-111 (2006).


\bibitem{hult1}{\sc A. Hultman,}
\newblock{\em Fixed points of involutive automorphisms of the Bruhat order,}
\newblock Adv. Math. 195, no. 1, 283--296 (2005).

\bibitem{kr}{\sc F.\ Knop, G. R\"ohrle,}
\newblock{\em Spherical subgroups in simple algebraic groups,}
\newblock arxiv:1305.3183


\bibitem {Kr}{\sc H.~Kraft,}
\newblock{\em  Parametrisierung von Konjugationklassen in ${\mathfrak sl}_n$},
\newblock {Math. Ann.}, {234} 209-220, (1978).

\bibitem{lu}{\sc J-H. Lu,}
\newblock{\em On a dimension formula for twisted spherical conjugacy classes in semisimple algebraic groups,}
\newblock Math. Z. Open Access DOI 10.1007/s00209-010-0776-4  (2010).

\bibitem{lus-indag}{\sc G. Lusztig,}
\newblock{\em A class of irreducible representations of the Weyl group,}
\newblock Proc. Kon. Nederl. Akad. Indag. Math. 82, 323-335 (1979).

\bibitem{lusztig-inventiones}{\sc G. Lusztig}
\newblock{\em Intersection cohomology complexes on a reductive group,}
\newblock Invent. Math. 75, 205--272 (1984).

\bibitem{lusztig1}{\sc G. Lusztig,}
\newblock{\em From conjugacy classes in the Weyl group to unipotent classes,}
\newblock Represent. Theory 15, 494-530 (2011).

\bibitem{lusztig2}{\sc G. Lusztig,}
\newblock{\em  From conjugacy classes in the Weyl group to unipotent classes, II,} 
\newblock Represent. Theory 16, 189Ð211 (2012).

\bibitem{lustrata}{\sc G. Lusztig,}
\newblock{\em On conjugacy classes in a reductive group,}
\newblock  arxiv:1305.7168v5

\bibitem{lusp}{\sc G. Lusztig, N. Spaltenstein,}
\newblock{\em Induced unipotent classes,}
\newblock J. London Math. Soc. (2), 19, 41--52 (1979).


\bibitem{mcninch-sommers}{\sc G. McNinch, E. Sommers}
\newblock{\em Component groups of unipotent centralizers in good characteristic,}
\newblock J. Algebra 270(1), 288--306 (2003).





\bibitem{kostant}{\sc Pasiencier, S., Wang, H.-C.,}
\newblock{\em Commutators in a semi-simple Lie group,}
\newblock Proc. Am. Math. Soc. 13(6), 907--913 (1962)



\bibitem{per-row}{\sc Perkins, S.; Rowley, P.,}
\newblock{\em Length of involutions in Coxeter groups,}
\newblock Journal of Lie Theory 14, 69--71 (2004).



\bibitem{sev}{\sc A. Sevostyanov}
\newblock{\em A proof of the De Concini-Kac-Procesi conjecture II. Strictly transversal slices to conjugacy classes in algebraic groups,}
\newblock arXiv:1403.4108v3, (2014).


\bibitem{sommers}{\sc E. Sommers}
\newblock{\em A generalization of the Bala-Carter theorem for nilpotent orbits,}
\newblock Internat. Math. Res. Notices 11, 539-562 (1998).

\bibitem{spalt}{\sc N. Spaltenstein,}
\newblock{\em Classes Unipotentes et sous-groupes de Borel,}
\newblock Springer-Verlag, Berlin (1982).



\bibitem{springer}{\sc T.A. Springer}
\newblock{\em Linear Algebraic Groups, Second Edition}
\newblock Progress in Mathematics 9, Birkh{\"a}user (1998).









\end{thebibliography}
\end{document}